\documentclass[10pt]{article}

\usepackage{amsfonts}

\newtheorem{theo}{Theorem}[section]
\newtheorem{lm}[theo]{Lemma}

\newtheorem{defin}[theo]{Definition}
\newtheorem{ex}{Example}[section]

\newcommand{\Reali}{\ensuremath{\mathbb{R}}}
\newcommand{\naturali}{\ensuremath{\mathbb{N}}}
\newcommand{\T}{{\mathbb T}}
\newcommand{\X}{{\mathbb X}}
\newcommand{\Xtilde}{{\widetilde{\mathbb X}}}
\newcommand{\G}{{\mathcal G}}
\newcommand{\OI}{{\mathcal O}}

\begin{document}
\title{Almost Vanishing Polynomials for Sets of  Limited Precision Points.}
\author{Claudia Fassino\thanks{Dipartimento di Matematica, Universit\`a di Genova, via Dodecaneso 35, Genova, Italy ({\tt fassino@dima.unige.it}).}}
\date{}
\maketitle
\begin{abstract}  
From the numerical point of view, given a set $\X\subset \Reali^n$ of $s$ points whose coordinates are known with only limited precision, each set $\Xtilde$  of $s$ points whose elements differ from those of  $\X$ of a quantity less than the data uncertainty can be considered equivalent to $\X$. We present an algorithm that, given  $\X$ and a tolerance $\varepsilon$ on the data error, computes a set $\G$ of polynomials such that each element of $\G$ ``almost vanishing'' at $\X$ and at all its equivalent sets $\Xtilde$. Even if $\G$ is not, in the general case,  a basis of the vanishing ideal $\mathcal I(\X)$, we show that, differently from the basis of $\mathcal I(\X)$ that can be  greatly influenced by the data uncertainty, $\G$ can determine a geometrical configuration simultaneously characterizing the set $\X$ and all its equivalent sets $\Xtilde$.
\end {abstract}
\section{Introduction}
Let $P=\Reali[x_1,\dots,x_n]$ be the polynomial ring in $n$ indeterminates over the reals and let $\X=\{p_1,\dots,p_s\}$ be  a finite set  of points of $\Reali^n$. 
 
 It is well known~\cite{KrRob00, KrRob05} that the vanishing ideal  $\mathcal I(\X)\subset P$ of all polynomials which vanish at $\X$ can be described by  a Gr\"obner basis~\cite{BM}, if a term ordering  is chosen,  or by a border basis, if an appropriate basis  of the quotient space $P/\mathcal I(\X)$ is given. 
 
 However, it is also well known that  small perturbations of the points of $\X$  can cause structural changes in the bases of $\mathcal I(\X)$~\cite{KrRob05, Ste04} as illustrated in the following example. 

\begin{ex}\label{first_ex}
{\rm  Let $\sigma$ be the DegLex term ordering  with~$x>y$. Given  the set of points $\X=\{(1,1), (3,2), (5.1,3)\}$, the  $\sigma$-Gr\"obner basis $GB$ of the vanishing ideal~$\mathcal I(\X)$ is given by: 
\begin{eqnarray*}
GB=
\left\{
\begin{array}{l}
y^2-20x+37y-18\\
 xy-43x+81 y -39\\
 x^2-90.1x+172.2y-83.1
\end{array}\right.
\end{eqnarray*}
The set  $GB$ is also the  border basis of $\mathcal I(\X)$, founded on the set  $\OI=\{1,\,y,\,x\}$ whose residue classes span $P/\mathcal I(\X)$.

A slightly perturbation of the point $(5.1,3)$ leads to a  new set of points  $\Xtilde =\{(1,1),(3,2), (5,3)\}$.  The $\sigma$-Gr\"obner basis  $\widetilde{GB}$ of the vanishing ideal $\mathcal I(\Xtilde)$  is completely different from $GB$:
\begin{eqnarray*}
\widetilde{GB}= \left\{
\begin{array}{l}
x-2y+1, \\
 y^3-6y^2+11y -6,
\end{array}\right.
\end{eqnarray*}
 Further, all the border bases of $\mathcal I(\Xtilde)$   also present a structural discontinuity, since the residue classes of the set  $\OI$ do not  span the space $P/\mathcal I(\Xtilde)$.\hfill{$\diamondsuit$}}
 \end{ex}
In the previous example the structural changes happen since the input points $\X$ are {\bf almost} aligned, while the slightly perturbed points  $\Xtilde$ are {\bf exactly} aligned.  This example also illustrates that, if we deal with a set $\X$ of points known with limited precision, the exact bases of $\mathcal I(\X)$ could not highlight some pleasant geometrical properties almost satisfied by the points $\X$. 

In this paper we present an algorithm that computes, given a set of points known with limited precision, a set of polynomials  allowing to recognize if such points almost lye on a  particularly simple geometrical configuration.

\medskip
Given a set $\X$ of  points whose coordinates are known with limited precision, each $p$ of $\X$ represents a ``cloud" of points: every point $\widetilde p$ which differs from $p$  by less than the data uncertainty  can be considered  computationally equivalent to $p$. Analogously,  an input set obtained from $\X$ replacing some $p$ by its perturbation $\widetilde p$ can be considered an {\bf admissible perturbation} computationally equivalent to $\X$.  It is then clear that the knowledge of $\X$ with limited precision,  combined with the structural  discontinuity of a basis, points out 
 that a significant characterization of $\mathcal I(\X)$ can be  a very tricky problem. In fact the structure of a basis can drastically change choosing different admissible input sets and moreover a blindly choice of a basis can hidden significant geometrical properties of $\X$. For this reason exact methods applied  to limited precision data can produce meaningless results. 

\medskip
The problem of the characterization of the vanishing ideal of a set of perturbed points has been studied by several authors from different points of view. 
In~\cite{Sauer07}, Sauer describes a method, suitable for numerical computations, which computes a small degree algebraic variety containing the input points.  In~\cite{HKPP}, Heldt~{\it et al.} present  an algorithm, based on the singular value decomposition of matrices, that computes, without using explicitly the estimation of the data error, a set of polynomials which assume particularly small values at the input points.

In~\cite{AFT08}, Abbott~{\it et al.} present  an algorithm that computes, explicitly using the tolerance on the data error, a monomial set $\OI$ which, in most cases, is a basis of  $P/\mathcal I(\X)$ and of  $P/\mathcal I(\Xtilde)$ for all the admissible perturbations $\Xtilde$ so that  the $\OI$-border basis of all the vanishing ideals $\mathcal I(\Xtilde)$ can be obtained. 

\medskip
Given a set $\X$ of limited precision points and a tolerance $\varepsilon$ on the data error, we focus our attention 
on the possibility of simultaneously characterizing, with a set of polynomials, the set $\X$  together with all its  admissible perturbations.    
To this aim, we present an algorithm that  computes an order ideal $\OI$ and a polynomial set $\G$, whose supports  are defined by $\OI$, having the following properties.
\begin{enumerate}
\item  The elements of $\G$ are almost vanishing, w.r.t.~the  norm of their coefficient vectors, at $\X$ and at each admissible perturbation $\Xtilde$.
\item For each admissible perturbation $\Xtilde$, the set $\{r(\Xtilde) | r \in \OI\}$ consists of independent vectors, up to the first order error analysis.
\item For each leading term $t$ of $g \in \G$  there could be an admissible perturbation $\widehat \X_g$ such that   $t(\widehat \X_g)$ depends on $\{r(\widehat \X_g) | r \in \OI\}$.
\end{enumerate}
Condition 3 implies that for each $g \in \G$ there could exist a polynomial $\widehat g$ with the same support of $g$ and similar coefficients which vanishes at $\widehat \X_ g$. If it is the case, the  algorithm determines a geometrical structure, given by $g$, almost satisfied by all the admissible perturbations of $\X$ and similar to a geometrical structure, given by $\widehat g$, exactly satisfied by the  admissible perturbation $\widehat \X_g$. 

As illustrated in the numerical examples in Section 6, it can happen that  there exists a single admissible perturbation $\widehat \X$ satisfying the previous property for all the polynomials $g \in \G$. In this case it  is very natural to consider $\widehat \X$  as a possible exact input set, that is the input in absence of  data error. Moreover, even if in  general $\G$ is not a basis of  $\mathcal I(\X)$, it is analogously natural to consider $\G$   as a common characterization of all the admissible perturbations of~$\X$. 

Finally, once again as we will show in Section 6, it can happen that $\widehat \X$ turns out to be the exact zero set of 
the polynomials of $\G$ so that, in this case, $\G$ is a Gr\"obner basis of $\mathcal I(\widehat \X)$.

\medskip
There is an evident open problem regarding the algorithm and its results: the existence and possibly the determination of the admissible perturbation $\widehat \X$. 
We will show, once again in Section 6, that there are cases when $\widehat \X$ does not exist. Then an open problem is to find conditions for the existence of $\widehat \X$ and, in case of existence, to determine it explicitly. 
The numerical tests suggest that in case of non existence, this can be due to two possible causes: 
the algorithm does not recognize a possible element of $\OI$ or it 
detects  some geometrical configurations, close to the points of  $\X$, which are incompatible with each other.
The study of such open problem will be the subject of our future work.
 
 \medskip
The paper improves and formalizes a 2005 Preprint of the author~\cite{Fas05} and it is organized as follows. Section 2 shows some basic concepts.  Section 3 contains the description of the algorithm. Section 4 and Section 5 describe
the  numerical properties of $\OI$ and $\G$, respectively. Finally,  Section 6 presents some examples illustrating the behaviour of our method.
\section{Preliminaries}
In order to  formalize the  idea of perturbed points, we recall the definitions of empirical point  and  of admissible perturbation~\cite{Ste04, AFT08}.

\begin{defin}\label{empirPts}
Let $p=(c_1,\dots, c_n)$ be a point of $\Reali^n$ and let $\varepsilon=(\varepsilon_1,\dots,
\varepsilon_n)$, with each $\varepsilon_i \in \Reali^+$, be the vector of the 
componentwise tolerances. An {\bf empirical point} $p^\varepsilon$
is the pair $(p,\varepsilon)$, where we call $p$ the {\bf specified value}
and $\varepsilon$ the {\bf tolerance}. 
A points $\widetilde p=(\widetilde c_1,\dots, \widetilde c_n) \in \Reali^n$ is called an {\bf admissible perturbation} of $p$  if $\widetilde c_i= c_i+e_i$, $|e_i| < \varepsilon_i$, $i=1,\dots,n$.
\end{defin}

Given a finite set $\X^\varepsilon$ of empirical points all sharing the same tolerance $\varepsilon$, we can formalize the concept of a set $\Xtilde$  ``equivalent" to $\X$ w.r.t. the data accuracy.

\begin{defin}
Let $\X^{\varepsilon} = \{p_1^{\varepsilon},\dots, p_s^{\varepsilon}\}$ 
be a set of empirical points with uniform tolerance~$\varepsilon$ and 
with $\X \subset \Reali^n$.  A set of points $\Xtilde=
\{\widetilde{p_1},\dots,\widetilde{p_s}\} \subset \Reali^n$ is called an 
{\bf admissible perturbation} of $\X$ if each $\widetilde{p}_i$ is an admissible perturbation of~$p_i$.
\end{defin}
 
Finally,  we recall (see~\cite{KrRob00,KrRob05}) some basic concepts related to the polynomial ring $P=\Reali[x_1,\dots,x_n]$.
\begin{defin}\label{eval}
Let $\X= \{p_1,\dots,p_s\}$ be a non-empty finite set of points of $\Reali^n$ 
and let $G=\{g_1,\dots,g_k\}$ be a non-empty finite set of polynomials.
\begin{itemize}
\item The $\Reali$-linear map ${\rm{eval}}_{\X} : P \rightarrow \Reali^s$ 
defined by ${\rm{eval}}_{\X}(f) = (f(p_1),\dots,f(p_s))$
is called the {\bf evaluation map} associated to $\X$.
For brevity, we write $f(\X)$ to mean ${\rm eval}_{\X}(f)$.
\item The {\bf evaluation matrix} (or {\bf vector} if $k=1$) of $G$ associated to $\X$, written
$M_{G}(\mathbb X)$ (or $g_1(\X)$), is defined as having entry $(i,j)$
equal to $g_j(p_i)$.
\end{itemize}
\end{defin} 

\begin{defin}
Let $\T^n$ be the monoid of power products of $P$ and let 
$\OI$ be a non-empty subset of $\T^n$.
\begin{itemize}
\item The set $\OI$ is called an {\bf order ideal} if $\OI= \overline{\OI}$, where  $\overline{\OI}$ is the set of all power products in $\T^n$ which divide some power product of $\OI$.
\item Given an order ideal $\OI$,  the {\bf corner set} of $\OI$ is the set 
$$\mathcal C[\OI] =\{t\in \T^n : t \notin \OI, x_i | t \Rightarrow t/x_i \in \OI, i=1\dots n \}$$
\end{itemize}
\end{defin}

Later on we suppose the reader familiar with the concepts of Gr\"obner basis and border basis of a vanishing ideal. Regarding these arguments, the reader is referred to the literature (see, among the others,~\cite{BM,KrRob00,KrRob05}). 

\section{The Numerical Algorithm }

Before processing a set $\X$ of limited precision points, it is possible to mitigate some negative effects of the data uncertainty, replacing with a single representative point the elements of $\X$ which differ each other by  less than the data accuracy, since they can be regarded  as different perturbations of the same empirical value. 
Later on we suppose w.l.o.g.~that the set $\X$ does not present such ``redundancy".  If it is not the case, it is possible to preprocess the input data to obtain well-separated points, using for instance the algorithms described in \cite{AFT07} and included in CoCoALib~\cite{Co}. 
Nevertheless the preprocessing of the input data is not sufficient to eliminate the instabilities of the exact bases of the vanishing ideal $\mathcal I(\X)$, as illustrated in Example~\ref{first_ex}, where the points $\X$ are well separated. 

\medskip
We base the construction of our algorithm  on the Buchberger-M\"oller  one~\cite{BM} which computes, given a set $\X$ of points and a  term ordering $\sigma$, the $\sigma$-Gr\"obner basis $GB$ of $\mathcal I(\X)$ as follows. At each step, if  $\OI=\{t_1,\dots,t_k\}$ is the order ideal computed at the previous steps, a power product  $t>_\sigma t_i$ is chosen. If  the vector $t({\X})$ is linearly independent of the vectors $\{t_1( {\X}),\dots,t_k(\X)\}$, $t$ is added to $\OI$.  Otherwise,  the  polynomial $g=t-\sum_{i=1}^k c_i t_i$ is put into $GB$. Nevertheless, since the test of linear dependence is crucially affected by even very small variations of the input data, when we deal with  points known with limited precision, small perturbations of the input data may lead to different choices in the Buchberger-M\"oller algorithm. 

In order to solve this drawback, we present an algorithm which  checks the linear dependence in a robust way w.r.t.~the data uncertainty. Since every admissible perturbation $\Xtilde$ is computationally equivalent to $\X$, the vector $t(\X)$ can be considered numerically dependent  on $\{t_1(\X),\dots,t_k(\X)\}$ if there exists an $\Xtilde$ such that  $t(\Xtilde)$ exactly depends on the vectors $\{t_1(\Xtilde),\dots,t_k(\Xtilde)\}$. Formally we have the following definition.
\begin{defin}\label{num_dep}
Given a set $\OI=\{t_1,\dots,t_k\}$ and a power product  $t$, the vector~$t(\X)$ {\bf numerically depends} on $\{t_1(\X),\dots,t_k(\X)\}$ if there exists an admissible perturbation $\Xtilde$ of $\X$ such that  the residual $\rho(\Xtilde)$ of the least squares problem $M_\OI(\Xtilde ) \widetilde \alpha = t(\Xtilde)$ is a null vector.
\end{defin}

\subsection{Sensitivity of the least squares problem}
In order to detect the numerical linear dependency of a set of evaluation vectors, we need some results concerning the sensitivity of the least squares problem
\begin{eqnarray}\label{main_problem}
M_\OI(\X) \alpha = t(\X)
\end{eqnarray}

First of all we recall some results, based on the componentwise perturbation analysis  (see~\cite{Hig94}), about the sensitivity of a generic least squares problem. 

Given an $h\times k$ matrix $A$,  we denote by $A^+$ its pseudoinverse, that is $A^+=(A^tA)^{-1}A^t$,
and by $|A|$ the matrix consisting of the absolute values of the elements of $A$; given an $h\times k$  matrix $B$, we assume that $|A| <|B|$ means that the relation holds componentwise. Moreover, given a real value $\eta \ll 1$ we denote by $O(\eta^m)$, $m\in \naturali$, an $h\times k$ matrix $W(\eta)=(w_{i,j}(\eta))$ (or a real function if $h=k=1$) such that, for each $(i,\; j)$,  $ |w_{i,j}(\eta)| / \eta^m$ is bounded near $0$.

\begin{theo}\label{generic}
Let  $A$ and $A+\Delta A$ both be $p\times q$, $p>q$, full rank matrices and let $b$ and $b+\Delta b$ be two vectors of $\Reali^p$ such that $|\Delta A | \le \eta E $ and $|\Delta b | \le \eta f$, $\eta\in\Reali^+$, $\eta \ll 1$.
Consider the least squares problems
\begin{eqnarray*}
 A x= b & with &  \rho=b-Ax  \;\; and \\
(A+\Delta A) (x+\Delta x)= b +\Delta b & with  & \rho+\Delta \rho=b+\Delta b -(A+\Delta A)(x+\Delta x)
\end{eqnarray*}
We have that
\begin{eqnarray*}
\Delta x &=& A^+ (\Delta b - \Delta A x ) + (A^tA)^{-1}(\Delta A)^t \rho + O(\eta^2)\\
\Delta \rho&=& (I- AA^+) (\Delta b - \Delta Ax) - ( A^+)^t (\Delta A)^t \rho + O(\eta^2)
\end{eqnarray*}  
where $I$ is the $p \times p$ identity matrix.
\end{theo}

\noindent{\bf Proof.}
It is possible~(see~\cite{GVL}) to express the least squares problem in the form
\begin{eqnarray*}
\left[\begin{array}{cc}
 I &  A \\
A^t  & 0 
\end{array}\right] \left[\begin{array}{c}
\rho\\
x
\end{array}\right] =\left[\begin{array}{c}
b \\
0 
\end{array}\right] 
\end{eqnarray*}
and so
\begin{eqnarray*}
\left[\begin{array}{cc}
 I & ( A+ \Delta A) \\
(A+ \Delta A)^t  & 0 
\end{array}\right] \left[\begin{array}{c}
\rho+\Delta \rho\\
x +\Delta x
\end{array}\right] =\left[\begin{array}{c}
b  + \Delta b\\
0 
\end{array}\right] 
\end{eqnarray*} 
 Taking the difference of the previous equations, we have 
\begin{eqnarray*}
\left[\begin{array}{cc}
 I &  A \\
A^t  & 0 
\end{array}\right] \left[\begin{array}{c}
\Delta \rho\\
\Delta x
\end{array}\right] =\left[\begin{array}{c}
\Delta b -\Delta A (x+\Delta x)\\
- (\Delta A)^t(\rho+\Delta \rho) 
\end{array}\right] \end{eqnarray*}  
Since
\begin{eqnarray*}
\left[\begin{array}{cc}
 I- AA^+ &  (A^+)^t \\
A^+  &  - (A^tA)^{-1}
\end{array}\right] 
{\rm is \; the \; inverse \; matrix \; of}
\left[\begin{array}{cc}
 I &  A \\
A^t  & 0 
\end{array}\right] 
\end{eqnarray*}  
we obtain 
\begin{eqnarray}
 \Delta \rho&=& ( I- AA^+) (\Delta b - \Delta A (x + \Delta x))- (A^+)^t (\Delta A)^t(\rho+\Delta \rho) \label{generale}\\
\Delta x &=&  A^+ (\Delta b - \Delta A (x + \Delta x)) + (A^tA)^{-1} (\Delta A)^t(\rho+\Delta \rho) \nonumber
 \end{eqnarray}  
Supposing  $|\Delta A | \le \eta E $ and $|\Delta b | \le \eta f$, the absolute values of $\Delta x$ and $\Delta \rho$ satisfy
\begin{eqnarray*}
|\Delta \rho | &\le& \eta \left ( \left | I- AA^+ \right | (f + E |x + \Delta x|) + | A^+|^t E^t|\rho+\Delta \rho| \right ) \\
|\Delta x | &\le& \eta \left (|A^+| (f + E |x + \Delta x|) + |(A^tA)^{-1}|E^t|\rho+\Delta\rho| \right ) 
 \end{eqnarray*}  
so we have that
$$ \Delta A \Delta x = O(\eta^2) \;\;\; {\rm and}  \;\;\; (\Delta A)^t \Delta \rho = O(\eta^2)$$
and the conclusion follows. \hfill$\diamondsuit$

\medskip
Since we are interested in the behaviour of the least squares problem~(\ref{main_problem}), we present an estimation of the  sensitivity of the matrix $M_\OI(\X)$ and of the vector $t(\X)$ to slight perturbations of the set~$\X$.

Given the power product $t=x_1^{\beta_1}\dots x_n^{\beta_n}$ and the monomial set  $\OI$, we denote by $\varepsilon_M=\max\{\varepsilon_i, \; i=1\dots n\}$, by $ \deg(x_k,t)=\beta_k$ the degree of $x_k$ into $t$, by
$\partial_k t = \deg(x_k,t) x_1^{\beta_1}\dots x_k^{\beta_k-1}\dots x_n^{\beta_n}$ and  by $  \partial_k \OI =  \left\{\partial_k t : t \in \OI\right\} $.

\smallskip
The following result concerns the sensitivity of the evaluation vector $t(\X)$.
\begin{lm}\label{lemma1}
Let $t$ be a power product of  $\T^n$.  Given a set $\X^\varepsilon=\{p_1^\varepsilon,\dots,p_s^\varepsilon\}$  of  empirical points and an admissible perturbation $\Xtilde=\{\widetilde p_1,\dots,\widetilde p_s \}$ of $\X$,  we have that the vector $\Delta t=t(\Xtilde) - t(\X)$ satisfies
\begin{eqnarray*}
\Delta t &=&\sum_{k=1}^n E_k  \partial_k t(\X)  +O(\varepsilon_M^2) 
\end{eqnarray*}
where  $E_k=Diag(e_{1,k},\dots, e_{s,k})$ is a diagonal matrix and $e_{i,k}$, $|e_{i,k} | <\varepsilon_k$, is a perturbation the $k$-th coordinate of $p_i$.
\end{lm}
\noindent{\bf Proof.}
First of all we consider  a  point $p=(c_1,\dots,c_n)\in \Reali^n$ and an admissible perturbation $\widetilde p=(\widetilde c_1,\dots,\widetilde c_n)$ of $p$ w.r.t the tolerance $\varepsilon$. Given $t=x_1^{\beta_1}\dots x_n^{\beta_n}$, we  have that
\begin{eqnarray*}
&&t(\widetilde p) -t(p) = ( c_1+ e_1)^{\beta_1}\dots ( c_n+e_n)^{\beta_n} -  c_1^{\beta_1}\dots c_n^{\beta_n}= \\
&&\sum_{k=1}^n e_k \beta_k c_k^{\beta_k-1} \prod_{h=1, h\neq k}^n c_h^{\beta_h} +O(\varepsilon_M^2) =
\sum_{k=1}^n e_k\partial_k t(p) +O(\varepsilon_M^2) 
\end{eqnarray*}
then we obtain
\begin{eqnarray*}
 t(\widetilde p_i) - t(p_i)&=&\sum_{k=1}^n e_{i,k} \partial_k t(p_i)  +O(\varepsilon_M^2)  
\end{eqnarray*}
and so, since $ t(\widetilde p_i) - t(p_i)$ is the $i$-th coordinate of $ t(\Xtilde) - t(\X)$, we conclude that
\begin{eqnarray*}
\qquad \qquad  \qquad \qquad t(\Xtilde) - t(\X)&=&\sum_{k=1}^n E_k  \partial_k t(\X)  +O(\varepsilon_M^2)  \qquad \qquad \qquad \qquad \diamondsuit
\end{eqnarray*}

The following result concerns the sensitivity of the evaluation matrix $M_\OI(\X)$.
\begin{lm}\label{lemma2}
Let $\OI$ be an order ideal.  Given a set $\X^\varepsilon=\{p_1\dots,p_s\}$  of  empirical points and an admissible perturbation $\Xtilde=\{\widetilde p_1,\dots,\widetilde p_s \}$ of $\X$,  we have that the matrix $\Delta M=M_\OI(\Xtilde) - M_\OI(\X)$ satisfies
\begin{eqnarray*}
\Delta M&=&\sum_{k=1}^n E_k  M_{\partial_k \OI} (\X)+O(\varepsilon_M^2) 
\end{eqnarray*}
where $E_k=Diag(e_{1,k},\dots, e_{s,k})$ is the diagonal matrix of Lemma~\ref{lemma1}. 
\end{lm}
\noindent{\bf Proof.}
Since the $j$-th column of $M_\OI(\X)$ is given by $t_j(\X)$, $t_j \in \OI$, Lemma~\ref{lemma1} implies that the $j$-th column of $M_\OI(\Xtilde) - M_\OI(\X)$ is 
\begin{eqnarray*}
t_j(\Xtilde) - t_j(\X)= \sum_{k=1}^n E_k  \partial_k t_j(\X)  +O(\varepsilon_M^2)
\end{eqnarray*}
The conclusion follows since $\partial_kt_j(\X) $ is  the $j$-th column of the evaluation matrix of the set  $\partial_k \OI$. \hfill$\diamondsuit$

\medskip
The next theorem, based on Theorem~\ref{generic}, Lemma~\ref{lemma1} and Lemma~\ref{lemma2}, presents a componentwise estimation of the sensitivity of the problem~(\ref{main_problem}) to the data perturbations.  Further, it
shows a componentwise upper bound of the absolute value of the residual, when  there exists  an admissible perturbation $\widehat \X$ such that the perturbed least squares problem  $M_\OI(\widehat\X) \widehat \alpha = t(\widehat\X)$ has a zero residual.

\begin{theo}\label{sens}
Let  $\X^\varepsilon$ be a set of $s$ empirical points and  let $\Xtilde$ be an admissible perturbation of $\X$.
Let $\OI$ be an order ideal such that $M_\OI(\X)$ and $M_\OI(\Xtilde)$ are full rank matrices.
Given the least squares problems
$$ M_\OI(\X) \alpha= t(\X) \;\; with \; residual \;\; \rho(\X)= t(\X) - M_\OI(\X) \alpha$$ 
and
$$M_\OI(\Xtilde) \widetilde \alpha= t( \Xtilde) \;\; with \; residual \;\; \rho(\Xtilde)= t(\Xtilde) - M_\OI(\Xtilde) \widetilde\alpha$$ 
then the vectors $\Delta \alpha= \widetilde \alpha - \alpha$ and $\Delta \rho=\rho(\Xtilde)-\rho(\X)$ satisfy
 \begin{eqnarray}\label{ro}
\begin{array}{rcl}
 \Delta \rho &= & (I- M_\OI(\X)M_\OI^+(\X)) \sum_{k=1}^n E_k \left( \partial_k t(\X)  -  M_{\partial_k \OI} (\X)  \alpha \right ) \\ \\
 &-& (M_\OI^+(\X))^t \left (\sum_{k=1}^n  M_{\partial_k \OI} ^t(\X) E_k\right ) \rho(\X)+O(\varepsilon_M^2)
 \end{array}
 \end{eqnarray}
\begin{eqnarray}\label{alfa}
 \begin{array}{rcl} 
\Delta \alpha &=&  M_\OI^+(\X) \sum_{k=1}^n E_k \left ( \partial_k t(\X)  -  M_{\partial_k \OI} (\X)  \alpha \right ) \\ \\
&+& (M_\OI^t(\X)M_\OI(\X))^{-1} \left ( \sum_{k=1}^n   M_{\partial_k \OI}^t (\X)  E_k \right )\rho(\X) + O(\varepsilon_M^2) 
\end{array}  
\end{eqnarray}  
Moreover,  if  there exists an admissible perturbation  $\widehat \X$ of $\X$ such that
the residual $\rho(\widehat \X)$ of the least squares problem  $M_\OI(\widehat \X) \widehat \alpha= t(\widehat \X)$ is a zero vector, then  the residual $\rho(\X)$  satisfies
\begin{eqnarray*}
\left | \rho(\X) \right | \le \left | I -M_\OI(\X)M_\OI^+(\X)  \right| \sum_{k=1}^n \varepsilon_k\left |\partial_k t(\X) - M_{\partial_k \OI}(\X)   \alpha \right | + O(\varepsilon_M^2)
\end{eqnarray*}
\end{theo}

\noindent{\bf Proof.} Since $\left | \Delta M \right| < \varepsilon_M E$, $\left | \Delta t \right| < \varepsilon_M f$ and $M_\OI(\X)$ and $M_\OI(\Xtilde)$ have full rank, from Theorem~\ref{generic} we obtain
\begin{eqnarray*}
\Delta \rho  &=& ( I- M_\OI(\X) M_\OI^+(\X)) (\Delta t - \Delta M \alpha )- ( M_\OI^+(\X))^t (\Delta M)^t\rho(\X)  + O(\varepsilon_M^2)\\
\Delta \alpha &=&  M_\OI^+(\X) (\Delta t - \Delta M  \alpha ) + (M_\OI^t(\X)M_\OI(\X))^{-1} (\Delta M)^t\rho(\X)  + O(\varepsilon_M^2)
\end{eqnarray*}  
and so, from Lemma~\ref{lemma1} and Lemma~\ref{lemma2}, Equations~(\ref{ro}) and (\ref{alfa})  follow.

Moreover, if $\rho(\widehat \X)$ is a zero vector from formula~(\ref{generale}) we have 
 \begin{eqnarray*}
 \rho(\X)=-\Delta \rho = ( M_\OI(\X)M_\OI^+(\X)-I) \sum_{k=1}^n E_k \left ( \partial_k t(\X)  -   M_{\partial_k \OI} (\X)  \alpha \right ) +O(\varepsilon_M^2) 
\end{eqnarray*}  
and if we consider the componentwise  absolute value of $ \rho(\X)$ we obtain
 \begin{eqnarray*}
\quad \left| \rho(\X) \right | &\le& \left | I- M_\OI(\X)M_\OI^+(\X)\right |\sum_{k=1}^n |E_k| \left | \partial_k t(\X)  -  M_{\partial_k \OI} (\X) \alpha  \right |  + O(\varepsilon_M^2)  \\
&\le& \left |I- M_\OI(\X)M_\OI^+(\X)\right |\sum_{k=1}^n \varepsilon_k  \left | \partial_k t(\X)  -  M_{\partial_k \OI} (\X)  \alpha \right |  + O(\varepsilon_M^2) \quad \diamondsuit
\end{eqnarray*}  

\subsection{The NBM Algorithm}
Theorem~\ref{sens} shows a sufficient condition for the numerical independency of $t(\X)$ of 
the columns of $M_\OI(\X)$. In fact if the residual $\rho(\X)$ of the least squares problem  $M_\OI(\X )\alpha = t(\X)$ satisfies
\begin{eqnarray}\label{cond}
 | \rho(\X) | >  \left |  I-M_\OI(\X)M_\OI^+(\X) \right | \sum_{k=1}^n \varepsilon_k \left|\partial_k t(\X) + M_{\partial_k \OI}(\X) \alpha\right | + O(\varepsilon_M^2)
 \end{eqnarray}
 then there are no admissible perturbations $\Xtilde$ of $\X$ such that the residual of the least squares problem $M_\OI(\Xtilde ) \widetilde \alpha = t(\Xtilde)$ is a null vector. So from Definition~\ref{num_dep} it follows that $t(\X)$ is numerically independent of~$\{ r(\X) : r \in \OI \}$. In particular, this implies that if $M_\OI(\Xtilde) $ is a full rank matrix then $[M_\OI(\Xtilde) t(\Xtilde)]$ is a full rank matrix too, for each admissible perturbation $\Xtilde$.
By exploiting this idea, we develop the {\bf N}umerical {\bf B}uchberger {\bf M}\"oller algorithm, whose main check is based on condition~(\ref{cond}). In particular, since we assume   the tolerance on the data error is relatively small, we  neglect the errors of order $O(\varepsilon_M^2)$ focusing our  attention  on  a first order error analysis of the problem.  

\vspace{3 mm}
\noindent{\bf The  Numerical Buchberger M\"oller (NBM) Algorithm.}

{\bf Input.} A set $\X^\varepsilon$ of $s$ empirical points and  a term ordering $\sigma$.

{\bf Output.} An order ideal  $\OI$   and a polynomial set $\G$.

\vspace{2 mm}
At the first step $\OI=\{1\}$ and $\G$ is an  empty set. A generic step can be described as follows. 
Let  $\OI=\{ t_1,\dots, t_k\} $ be the order ideal computed at the previous steps
and let $t$ be the current power product, $t>_\sigma t_1,\dots, t_k $, chosen according to the strategy of 
the Buchberger-M\"oller algorithm.
\begin{enumerate}
 \item Solve the least square problem $M_\OI(\X)\alpha=t(\X)$ and compute the residual   $\rho(\X)=t(\X)-M_\OI(\X) \alpha$.
 \item  If $\rho(\X)$ satisfies
\begin{eqnarray}\label{condition}
| \rho(\X) |  >\left | I - M_\OI(\X) M_\OI^+(\X) \right | \sum_{k=1}^n \varepsilon_k \left |\partial_k t(\X) - M_{\partial_k \OI}(\X)  \alpha \right | 
\end{eqnarray}
 then put  the term $ t$ into the set $\OI$.
 \item Otherwise, put the polynomial  $g=t - \sum_{i=1}^k \alpha_i t_i$ into $\G$.\hfill $\diamondsuit$
\end{enumerate}

The NBM algorithm stops after finitely many steps and computes  an order ideal $\OI$, since the strategy to choose the power products to analyze  is the same as in the Buchberger-M\"oller algorithm. 

Note that the term ordering $\sigma$ is only a computational tool for obtaining a set $\OI$ closed under taking divisors.
In fact in the general case $\OI$ is different from $\OI_\sigma$,  the quotient basis determined by the $\sigma$-Gr\"obner basis of $\mathcal I(\X)$. Moreover it can happen that, for each possible term ordering $\tau$, $\OI$ does not coincide  to any $\OI_\tau$ corresponding to the $\tau$-Gr\"obner basis of $\mathcal I(\X)$ (see Example~\ref{ex4}). For this reason any different strategy for building an order ideal can be used in the NBM algorithm instead to fix a term ordering. 
\section{Properties of the order ideal $\OI$}
A first important property of  the order ideal $\OI$ computed by the NBM algorithm is its invariance w.r.t.~the scaling and the translation of the points $\X$, as shown in the following theorem.
\begin{theo}\label{scaling}
Let  $\X^\varepsilon$  be a set of empirical points with
\begin{eqnarray*}
\X=\{p_1,\dots,p_s\}\;\;\; p_i=(c_{i,1},\dots,c_{i,n}) \;\;\; and \;\;\;\varepsilon= (\varepsilon_1,\dots,\varepsilon_n) 
\end{eqnarray*}
Let  $\X_S^{\delta}$  be the set of scaled empirical points such that $\X_S=\{\overline p_1,\dots, \overline p_s\}$ and
\begin{eqnarray*}
\overline p_i=(d_1 c_{i,1},\dots,d_n c_{i,n}) \;\;\;  with \;\;\; (d_1,\dots,d_n) \in \Reali^n \;\;\; and \;\;\; \delta= (|d_1| \varepsilon_1,\dots,|d_n|\varepsilon_n)
\end{eqnarray*}
Let  $\X_T^{\tau}$  be the set of translated empirical points such that $\X_T=\{\widehat p_1,\dots, \widehat p_s\}$ and
\begin{eqnarray*}
\widehat p_i=(c_{i,1}+v_1,\dots, c_{i,n}+v_n) \;\;\;  with \;\;\; (v_1,\dots,v_n)  \in \Reali^n \;\;\; and \;\;\; \tau= \varepsilon
\end{eqnarray*}
Then the NBM algorithm computes the same order ideal $\OI$ for all the input sets $\X^\varepsilon$, $\X_S^\delta$ and $\X_T^\tau$.
\end{theo}
\noindent{\bf Proof.}
We prove that the NBM algorithm computes the same order ideals at each step independently of the input sets
$\X^\varepsilon$, $\X_S^{\delta}$ or $\X_T^\tau$.

 At the first step it is true, since $\OI=\{1\}$.
Let us suppose that, at the current step, with  all the three  input sets the same order ideal $\OI=\{t_1,\dots,t_k\}$ has been computed and that the term $t$ has to be processed.

\smallskip
 Let us consider the set $\X_S$ of the scaled points. 
 
 Given a term  $r=x_1^{\beta_1}\dots x_n^{\beta_n}$, denoting by $r(d)=d_1^{\beta_1}\dots d_n^{\beta_n}$, we have
\begin{eqnarray*}
 r(\overline p_i) = r(d) r(p_i)\;\;\;\; {\rm and} \;\;\;\;  \partial_k r(\overline p_i)= \frac{r(d)}{d_k}\partial_k  r( p_i)
 \end{eqnarray*}
 so that,  denoting by $D_\OI$ the diagonal matrix $Diag(r(d) : r \in \OI)$,
\begin{eqnarray*}
 t(\X_S)=t(d)t(\X) &{\rm and}& M_\OI(\X_S) = M_\OI(\X) D_\OI \\
 \partial_kt(\X_S)=\frac{t(d)}{d_k}\partial_k t(\X) &{\rm and} &M_{\partial_k \OI}(\X_S) = \frac{1}{d_k}M_{\partial_k\OI}(\X) D_\OI 
 \end{eqnarray*}
The least squares problems $M_\OI( \X)  \alpha= t(\X)$ and $M_\OI( \X_S) \alpha_S= t(\X_S)$ solved with input sets $\X^\varepsilon$ and $\X_S^\delta$ are such that 
\begin{eqnarray*}
&M_\OI( \X)  D_\OI \alpha_S =  t(d) t(\X) \;\; \Rightarrow \;\;   D_\OI \alpha_S =  t(d) M^+_\OI( \X) t(\X)   \;\; \Rightarrow \;\;  \alpha_S = t(d) D^{-1}_\OI \alpha &\\
&\rho(\X_S) = t( \X_S) - M_\OI(\X_S) \alpha_S = t(d)t(\X)- t(d)M_\OI(\X)D_\OI D_\OI^{-1} \alpha = t(d) \rho(\X)&
 \end{eqnarray*}
If we consider the upper bound~(\ref{condition}) of Step 2  computed for the scaled empirical points $\X_S^\delta$, straightforward computations show that
 \begin{eqnarray*}
I-M_\OI(\X_S)M_\OI^+(\X_S) & = &I - M_\OI(\X)M_\OI^+(\X) \\
\partial_k t(\X_S) - M_{\partial \OI_k}(\X_S) \alpha_S &=& \frac{t(d)}{d_k} 
\Big [\partial_k t(\X) - M_{\partial \OI_k}(\X) \alpha \Big ] 
 \end{eqnarray*}
It follows that $t$ satisfies condition~(\ref{condition}) with input set $\X_S^\delta $  if and only if 
 \begin{eqnarray*}
|t(d)| | \rho( \X) |  >|t(d)| \left | I - M_\OI(\X) M_\OI^+(\X) \right | \sum_{k=1}^n \frac{\delta_k}{|d_k|} \Big |\partial_k t(\X) - M_{\partial_k \OI}( \X)  \alpha \Big | 
\end{eqnarray*}
that is if and only if $t$ satisfies condition~(\ref{condition}) with input set $\X^\varepsilon$ since $\delta_k = |d_k| \varepsilon_k$. We conclude that the NBM algorithm puts $t$ into $\OI$ processing  the input $\X^\varepsilon$ if and only if $t$ is added to $\OI$ using the input  $\X_S^\delta$.

\smallskip
 Let us consider the set $\X_T$ of the translated points.  
 
 Given a term  $r=x_1^{\beta_1}\dots x_n^{\beta_n}$,  there exist~(see~\cite{Tor08}) a set $R=\{r_j : r_j | r\}$ of  power products  and a set $\{ \gamma_j : \gamma_j=\gamma_j(v_1,\dots,v_n)\}$ of coefficients such that for each $p=(c_1,\dots,c_n)$  and $\widehat p=(c_1+v_1,\dots,c_n+v_n)$ 
 \begin{eqnarray*}
 r(\widehat p) = r(p) +\sum_{r_j \in R} \gamma_j r_j(p) 
\end{eqnarray*}
Furthermore, let $F_{(v_1,\dots,v_n)}: \Reali^{n}\rightarrow \Reali$ be a function  such that 
$$F_{(v_1,\dots,v_n)}(x_1,\dots,x_n) =( x_1+v_1)^{\beta_1}\dots ( x_n+v_n)^{\beta_n} -  x_1^{\beta_1}\dots x_n^{\beta_n} - \sum_{r_j \in R} \gamma_j r_j(x_1,\dots,x_n) $$
Since $F_{(v_1,\dots,v_n)}(p)=0$ for each point $p\in \Reali^n$ we obtain
$$0=\frac{\partial F_{(v_1,\dots,v_n)}}{\partial x_k}(p) =\frac{\partial r}{\partial x_k} (\widehat p) -\frac{\partial r}{\partial x_k}(p)  -  \sum_{r_j \in R} \gamma_j\frac{\partial r_j}{\partial x_k}(p) $$
that is, using our notation, 
$$ \partial_k r (\widehat p)  = \partial_kr(p)  +  \sum_{r_j \in R} \gamma_j \partial_k r_j(p) $$
Now, let us consider at the current step the set $\OI$ and the power product $t$. By  construction $\OI$ is factor closed,   so that for each $ r \in \OI \cup \mathcal C[\OI]$ the set $R$ is a subset of $\OI$. 
Since $t \in \mathcal C[\OI]$,  we have
\begin{eqnarray*}
t(\widehat p_i)= t(p_i) + \sum_{t_j \in \OI} \lambda_j t_j(p_i)\;\;\; {\rm and}\;\;\; \partial_k t(\widehat p_i)= \partial_k t(p_i) + \sum_{t_j \in \OI} \lambda_j \partial_k t_j(p_i)
 \end{eqnarray*}
so that, denoting by $\lambda$ the vector which consists of the values $\lambda_j$,
\begin{eqnarray*}
t( \X_T)= t(\X) + M_\OI(\X) \lambda \;\;\; {\rm and}\;\;\; \partial_k t(\X_T)= \partial_k t(\X) + M_{\partial_k \OI}(\X)\lambda
 \end{eqnarray*}
Analogously, analyzing each column of the matrices $M_\OI(\X_T)$ and  $M_{\partial_k}\OI(\X_T)$, there exists a square matrix $\Lambda$ such that
\begin{eqnarray}\label{matrices}
M_\OI( \X_T)= M_\OI(\X) + M_\OI(\X) \Lambda \;\;\;{\rm and}\;\;\; M_{\partial_k \OI}(\X_T)= M_{\partial_k \OI}(\X) + M_{\partial_k \OI}(\X)\Lambda
 \end{eqnarray}
The least squares problems $M_\OI( \X)  \alpha= t(\X)$ and $M_\OI( \X_T)  \alpha_T= t( \X_T)$ solved with the  input sets $\X^\varepsilon$ and $\X_T^\tau$ are such that 
\begin{eqnarray*}
&& M_\OI(\X)(I+ \Lambda)\alpha_T  = t(\X)+M_\OI(\X) \lambda \;\; \Rightarrow \;\; (I+ \Lambda)\alpha_T = \alpha +\lambda \\
&&  \rho(\X_T)  = t( \X_T) - M_\OI(\X_T)  \alpha_T = t(\X) +M_\OI(\X)\lambda  - M_\OI(\X)(\alpha +\lambda )= \rho(\X)
 \end{eqnarray*}
 Since the residual of least squares problem is invariant w.r.t.~the translation and $M_\OI(\X)$ has full rank then 
  $M_\OI(\X_T)$ is a full rank matrix too. It follows from~(\ref{matrices}) that  $M_\OI(\X)(I+\Lambda)= M_\OI(\X_T)$, and so  $I+\Lambda$ is a non singular matrix.

If we consider the upper bound~(\ref{condition}) of Step 2,  computed for the translated  empirical points $\X_T^\tau$ straightforward calculations lead to
$$I-M_\OI(\X_T)M_\OI^+(\X_T) = I - M_\OI(\X)M_\OI^+(\X)$$
Furthermore since 
 \begin{eqnarray*}
&& \partial_k t(\X_T) - M_{\partial \OI_k}(\X_T)  \alpha_T =\partial_k t(\X)  + M_{\partial_k\OI}(\X) \lambda  - M_{\partial \OI_k}( \X)(I+ \Lambda)  \alpha_T =\\ 
&& \partial_k t(\X)  + M_{\partial_k\OI}(\X) \lambda  - M_{\partial \OI_k}( \X)(\alpha + \lambda) =   \partial_k t(\X)   - M_{\partial \OI_k}( \X)\alpha  
 \end{eqnarray*}
it follows that $t$ satisfies condition~(\ref{condition}) with input $ \X_T^\tau$  if and only if 
 $t$ satisfies condition~(\ref{condition}) with input  $\X^\varepsilon$. We conclude that the NBM algorithm puts $t$ into $\OI$ processing  $\X^\varepsilon$ if and only if $t$ is added to $\OI$ using~$\X_T^\tau$. \hfill$\diamondsuit$

\medskip
In order to analyze the stability properties of the order ideal $\OI$ we recall some basic concepts (see~\cite{AFT08}).

\begin{defin}
An order ideal~$\OI$ is {\bf stable} w.r.t.~$\X^\varepsilon$ if the evaluation matrix~$M_{\OI}(\Xtilde)$ has full
rank for each admissible perturbation~$\Xtilde$ of~$\X^\varepsilon$.
 \end{defin}

Heuristically speaking an order ideal $\OI$ can be considered stable w.r.t.~the data uncertainty if the linear independency  of the evaluation vectors of its elements is not affected by slight perturbations of  $\X$.

It is well known that each order ideal  is  stable  providing  the values of $\varepsilon$ are sufficiently small.
Nevertheless in our problem the tolerance $\varepsilon$ is given a priori  and then not all the order ideals turn out to be stable. 

\smallskip
By the very nature of the NBM algorithm, no formal results about the stability of $\OI$ can be stated. In fact, when the 
numerical independence of  $\{ r(\X) : r \in \OI\}$ is tested using  condition~(\ref{cond}), Theorem~\ref{sens} ensures that $\OI$ is stable. Unfortunately, for implementative reasons, the NBM algorithm checks the numerical independence of  $\{ r(\X) : r \in \OI\}$ using the first order approximation~(\ref{condition}) of~(\ref{cond}). So the stability of $\OI$ is not guaranteed.

Nevertheless, in most cases, the numerical tests show that the upper bound~(\ref{condition}) is satisfied with a wide margin, widely greater than $O(\varepsilon_M^2)$: then, although we cannot have the complete certainty, there is an high probability that the order ideal $\OI$ is stable. 

\smallskip
We recall that~(see~\cite{AFT08})  it  is possible to compute a stable order ideal by using the SOI algorithm.
Its elevated computational cost, widely greater than the computational cost of the NMB algorithm,  makes the SOI algorithm not particularly suitable for all who are not mainly interested in stability. 

However, the possibility of comparing  the results of the SOI and the NBM algorithms points out a comforting behaviour 
of the NBM algorithm. In fact  in several numerical tests the order ideals computed by the algorithms coincide. This result supports the fact  that  the NBM algorithm, although without certainty,  often gives stable order ideals.
\section{Properties of the polynomial set $\G$}
First of all, we formalize the idea of almost vanishing polynomials introducing the following definition. 
Theorem~\ref{scaling} allows to restrict our attention to set of points whose coordinate belong to $[-1,\;1]$.
\begin{defin}
Given a set $\X^\varepsilon$ of empirical points whose coordinates belong to $[-1,1]$ a polynomial $g$, with coefficient vector $c$, is {\bf almost vanishing} at $\X$ if
$$\frac{\|g(\X) \|_2}{\|c\|_2} < O(\varepsilon_M) $$
\end{defin}

Obviously in the general case $\G$ is not a basis of $\mathcal I(\X)$, since $\G$ can contain polynomials that  do not exactly vanish at $\X$. However the following theorem shows that $\G$ exhibits interesting properties w.r.t.~the data uncertainty. 

\begin{theo}\label{polynomial_set}
Let  $\X^\varepsilon$ be a set of $s$ empirical points and let $\Xtilde$ be an admissible perturbation of $\X$. The polynomial set   $\G$ satisfies the following properties.
\begin{description}
\item{\bf P1} If $g$ is a polynomial of $\G$ of degree $\deg(g)$  and coefficient vector $c$, then
\begin{eqnarray*}
\frac{ \| g(\X)\|_2}{\|c\|_2} < s \;\deg(g) \sum_{k=1}^n \varepsilon_k  \;\; and \;\; 
\frac{ \| g(\Xtilde)\|_2}{\|c \|_2} < 2 \; s \; \deg(g) \sum_{k=1}^n \varepsilon_k +O(\varepsilon_M^2)
\end{eqnarray*}
Therefore,  $g$ is  almost vanishing at $\X$ and at $\Xtilde$.
\item{\bf P2} If the zero set of $\G$  is an admissible perturbation $\widehat \X$ such that $M_\OI(\widehat \X)$ has full rank, then $\G$ is the $\sigma$-Gr\"obner basis of  $\mathcal I(\widehat \X)$.
\item{\bf P3} If $\#\OI=s$,  each polynomial $g\in\G$  corresponds to  a unique polynomial $ g_b$ of the  $\OI$-border basis of  $\mathcal I(\X)$ such that the support of $g$ is a subset of the support of $g_b$. Furthermore, if $c$ and $c_b$ are respectively  the coefficient vectors of $g$ and $g_b$ then
\begin{eqnarray*}
\frac{\| c_b- [c, 0\dots 0]\|_2}{\|c\|_2} \le \deg(g) \|M_\OI(\X)\|_2 \|M_\OI^{-1}(\X)\|_2 \sum_{k=1}^n \varepsilon_k \end{eqnarray*}
 \end{description}
\end{theo}

\noindent{\bf Proof.}
\begin{description}
\item{\bf P1} Let us consider Step 3 of the NBM algorithm where the polynomial $g$ is computed. Let $t$ be the monomial analyzed at such step and let $\OI_t$ be  the order ideal obtained at the previous ones. Since the polynomial $g$ is added to $\G$ if the residual $\rho(\X)$ of the least squares problem $M_{\OI_t}(\X) \alpha = t(\X)$ does not satisfy condition~(\ref{condition}) and since $g(\X)=\rho(\X)$ we have 
$$\| g(\X)\|_2 < \| I - M_{\OI_t}(\X) M_{\OI_t}^+(\X)\|_2 \sum_{k=1}^n \varepsilon_k \left \|\partial_k t(\X) - M_{\partial_k \OI_t}(\X)  \alpha \right \|_2$$
First of all we prove that  $\| I - M_{\OI_t}(\X) M_{\OI_t}^+(\X)\|_2=1 $. In fact let $A$ be a $p \times q $ full rank matrix $A$, $p > q$  and  let $A=U \Sigma V^t$ be its singular values decomposition. It is well known~(see~\cite{GVL}) that $U$ and $V$ are square orthonormal matrices and $\Sigma$ is a block matrix of the form $\Sigma^t=[\Sigma_1, \; 0]$, where $\Sigma_1$ is the square diagonal matrix of the singular values, and so $ \| I - A A^+ \|_2= \| I - \Sigma (\Sigma ^t \Sigma)^{-1} \Sigma \|_2 =1 $.

Later on we denote by $\widehat M_k$ the matrix $ [ \partial_k t(\X)  \; M_{\partial_k \OI_t} (\X)]$, which consists of the vectors $\partial_kr(\X)$ with $r \in \OI_t \cup \{t\}$.
Obviously, if $\deg(x_k,r)=0$ the corresponding column of $\widehat M_k$ is a null vector. Moreover, 
for each $q, r\in \OI_t \cup \{t\}$ such that $q\neq r$, $\deg(x_k,r)\neq 0$ and $\deg(x_k,q)\neq 0$, we have
$ {\partial_k r}/{\deg(x_k,r)} \neq {\partial_k  q}/{\deg (x_k, q)}$ and,  since $\OI_t$ is factor closed and $t \in \mathcal C[\OI]$, 
$ {\partial_k r}/{\partial (x_k,r)} \in \OI_t $.

It follows that each column $\partial_k r(\X)$ of $\widehat M_k$ is a null vector or it corresponds to a unique column   of $M_{\OI_t}(\X)$ multiplied by $\deg (x_k,r)$.
Since $\|\widehat M_k\|_2 $ is equal to the norm of its submatrix consisting of the non zero columns and since $\deg(x_k,r) \le \deg(g)$ we have that 
 $$\|\widehat M_k\|_2 \le  \deg(g)\|M_{\OI_t}(\X)\|_2$$
Finally, since $c=[1,\; -\alpha]^t$ is the coefficient vector of $g$, we have 
\begin{eqnarray}\label{emmecap}
\left\|\partial_k t(\X) - M_{\partial_k \OI_t} (\X) \alpha \right \|_2= \left\|  \widehat M_k c\right \|_2\le  
 \deg(g)\|M_{\OI_t}(\X)\|_2 \|c\|_2\end{eqnarray}
 so that 
 \begin{eqnarray}\label{magg}
\|g(\X)\|_2 \le    \deg(g)\|M_{\OI_t}(\X)\|_2 \|c\|_2 \sum_{k=1}^n \varepsilon_k 
 \end{eqnarray}
 Since the coordinates of the points belong to $[-1,\;1]$, we have $\|M_{\OI_t}(\X)\|_2 \le \|M_{\OI_t}(\X)\|_F \le s$, where $\| \cdot \|_F$ is the Frobenius matrix norm. Then the first upper bound of P1 follows immediately.

Further, in order to show the result about $g(\Xtilde)$, note that 
\begin{eqnarray*}
 g(\Xtilde)&=&g(\X) +  \Delta t - \Delta M\alpha
 \end{eqnarray*}
So from Lemma~\ref{lemma1} and Lemma~\ref{lemma2} we have
\begin{eqnarray*}
g(\Xtilde)&= &g(\X)+ \sum_{k=1}^n E_k \left( \partial_k t(\X)  -  M_{\partial_k \OI_t} (\X) \alpha \right )+O(\varepsilon_M^2)
\end{eqnarray*}
and, computing the norm of $g(\Xtilde)$,  we obtain 
\begin{eqnarray*}
\|g(\Xtilde)\|_2 &\le& \|g(\X)\|_2+\deg(g)   \|M_{\OI_t}(\X)\|_2\|c\| \sum_{k=1}^n \varepsilon_k  +O(\varepsilon_M^2) \\
& \le& 2\; s \; \deg(g)\|c\|_2 \sum_{k=1}^n \varepsilon_k +O(\varepsilon_M^2)
\end{eqnarray*}
that is the second upper bound of P1.
\item{\bf P2}  If the zero set of $\G$  is an admissible perturbation $\widehat \X$,  since the residuals associated to the elements of $\G$ vanish at $\widehat \X$ and $M_\OI(\widehat \X)$ has full rank, the NBM algorithm computes the polynomial set $\G$ with input set $\widehat \X$ and tolerance $\varepsilon=(0,\dots,0)$.
Then Property P2 follows immediately because the NBM algorithm with a zero tolerance coincides with the exact Buchberger-M\"oller one.
\item{\bf P3} Since $\# \OI=s$ then $\OI$ is the quotient basis of $P/\mathcal I(\X)$ and so there exists the $\OI$-border basis of $\mathcal I(\X)$~(see~\cite{KrRob05}). By construction, each polynomial $g\in \G$ with leading term $t$ and support contained in $\{t\} \cup  \OI_t$ corresponds to a polynomial $g_b $ of the $\OI$-border basis of $\mathcal I (\X)$  whose support is contained in $\{t\} \cup \OI$.
If we order  the elements of $\OI_t$ and $\OI$ in an increasing way w.r.t.~$\sigma$, then the columns of $M_{\OI_t}(\X)$ coincide with the first $\#\OI_t$ columns of $M_\OI(\X)$ and the coefficient vectors $c=[1, \;-\alpha]$ of $g$ and $c_g=[1,\;- \beta]$ of $g_b$ obey  $\|c_g-[c, 0 \dots,0]\|=\|\beta - [\alpha, 0 \dots 0]\|$. Moreover they are such that
\begin{eqnarray*}
M_{\OI}(\X)\left [\begin{array}{c} \alpha \\ 0 \end{array} \right ]= t(\X)+\rho(\X) \qquad \qquad \qquad M_{\OI}(\X) \beta = t(\X)
\end{eqnarray*}
Then we obtain 
\begin{eqnarray*}
\left [\begin{array}{c} \alpha \\ 0 \end{array} \right ]- \beta= M_{\OI}^{-1}(\X)\rho(\X) 
\end{eqnarray*}
From  $\rho(\X)=g(\X)$ and the upper bound~(\ref{magg}), the thesis follows.\hfill$ \diamondsuit$
\end{description}

\vspace{3 mm}
Let us point out two pleasant properties of the polynomial set $\G$ easily following by P1 and P3.

Property P1 immediately implies that $\G$ can contain 
almost vanishing polynomials even if the coordinates of $\X$ do not belong to $[-1,\; 1]$. In fact it is sufficient that  the coordinates of $\X$  are ``not too elevated'' (see the examples of Section 6).

In the case when  the condition number $ \|M_{\OI}^{-1}(\X)\|_2\|M_\OI (\X) \|_2 $~(see~\cite{GVL}) of the matrix $M_\OI(\X)$ is ``not too elevated'', Property P3 implies that  $g$ is  ``close" to a polynomial $g_b$ vanishing at $\X$.  Then $\X$  can be considered a {\bf  pseudozero set} of $\G$, in the sense given by Stetter~(see~\cite{Ste199, Ste299}).

The formal results proved above allow us to justify in details the heuristic properties of $\G$ described in the Introduction and in particular the reasons why $\G$  characterizes the input points $\X$.

First of all note that Property P1 implies that each element $g$ of $\G$ assumes small values, and then it is almost vanishing, at $\X$ and at each admissible perturbation (of course w.r.t.~the norm of its coefficient vector).

Moreover, we recall that, by construction, each element $g$ of $\G$ with leading term $t$ corresponds to a least squares problem $M_{\OI_t}(\X)\alpha = t(\X)$, $\OI_t \subset \OI$, whose residual $\rho$ does not satisfy  condition~(\ref{condition}). Note that, since Theorem~\ref{sens} involves only sufficient conditions on the residual, the fact that $\rho$ does not satisfy condition~(\ref{condition}) gives essentially the same information of the fact that $\rho$ does not satisfy condition~(\ref{cond}).

Given $g$ of $\G$,  let us suppose that there exists an admissible perturbation $\widehat \X_g$ such that $\rho(\widehat \X_g)$ is a null vector. This is a possible case because condition~(\ref{cond}) is not satisfied.  Moreover,  let us suppose that the order ideal $\OI$ is stable, so that the matrix $M_{\OI_t}(\widehat \X_g)$ has full rank. Then there exists a polynomial $\widehat g$, given by the solution of $M_{\OI_t}(\widehat \X_g)\widehat \alpha = t(\widehat\X_g)$, having the following properties:
\begin{itemize}
\item $\widehat g$ exactly vanishes at $\widehat \X_g$;
\item $\widehat g$ has the same support of $g$;
\item $\widehat g$ and  $g$ have ``similar'' coefficients, if the condition number $\|M_\OI(\X)\|\|M^+_\OI(\X)\|$ is ``not too elevated''. In fact from relations~(\ref{alfa}) and~(\ref{emmecap}) we have
\begin{eqnarray*}
 && \frac{\|[1,\; -\widehat \alpha] - [1,\; -\alpha]  \|_2}{\|[1,\; -\alpha]  \|_2} \le \frac{\sum_{k=1}^n \varepsilon_k \|\partial_kt(\X) - M_{\partial_k \OI} (\X) \alpha \|}{\|[1,\; -\alpha]  \|_2} \le  \\
 && \|M^+_\OI(\X)\|\|M_\OI(\X)\| \deg(g) \sum_{k=1}^n \varepsilon_k \end{eqnarray*}
\end{itemize}
In this sense $g$ can selects a geometrical configuration $\widehat \X_g$ of points, close to $\X$,
that can be considered an ``approximate" representation of the input points independent of the data errors.

Furthermore, as we will show in the examples of Section 6, it can happen that the whole set of polynomials $g$ of $\G$ selects a unique geometrical configuration~$\widehat \X$. Therefore the  polynomials $\widehat g$ constitute a Gr\"obner basis of $\mathcal I(\widehat \X)$. We can then conclude that $\G$ can be viewed as an approximation of a Gr\"obner basis of an ideal of points $\widehat \X$ close to $\X$ and the set $\widehat \X$ can be considered as a possible ``exact" configuration corresponding to the absence of data uncertainty.

 We point out that, once again as shown in Section 6, it can happen that each $\widehat g$ coincides with $g$ and then $\G$ itself is a Gr\"obner basis for $\mathcal I(\widehat \X)$.

\medskip
Let us conclude this section with a short recall of the open problems already presented in the Introduction.
They are essentially related to the existence and possibly the determination of $\widehat \X$.
The numerical examples show that there are cases when $\widehat \X$ does not exist. This seems to be due to two possible causes. One is because the NBM algorithm could not recognize a possible element of $\OI$ so that a polynomial $g$ which never vanishes at any admissible perturbation is added to $\G$. The second reason is when
 the points of  $\X$ are close to different incompatible  geometrical configurations. However, in our numerical examples,  in this case the NBM algorithm explicitly detects these incompatible geometrical configurations.
\section{Numerical examples}
The following numerical tests are performed using a prototype version of the NBM algorithm. An improved version of it will be included soon in CoCoALib~\cite{Co}.
In all the examples the term ordering DegLex, $y<x$ is used; in addition, the coordinates of the points  and the
coefficients of the polynomials  are displayed with a finite number of digits, but all computations are performed in exact
arithmetic using CoCoA~4.7~\cite{Co}. 

\medskip
In Example~\ref{ex1} the NBM algorithm computes an exact Gr\"obner basis of a vanishing ideal of an admissible perturbation.

 \begin{ex}\label{ex1}
{\rm Given the same data of  Example~\ref{first_ex}, that is  the set of  points $\X=\{(1,1), (3,2), (5.1,3)\}$, if the tolerance is  $\varepsilon=(0.15,\;0)$, the NBM algorithm  computes the quotient basis $\OI=\{1,y,y^2\}$ and the polynomial set $\G$: 
\begin{eqnarray*}
\G = \left\{
\begin{array}{rcl}
g_1&=&x-2.05y+1.0{\overline 6} \\
g_2&=&y^3-6y^2+11y-6
\end{array}
\right. 
\end{eqnarray*}
The polynomial  $g_2$ vanishes at $\X$ while   $g_1$, with   coefficient vector $c_1$,  is almost vanishing at $\X$ since $\|g_1(\X)\|/\|c_1\|=0.0162$. 

Since  $\G$ is the  $\sigma$-Gr\"obner basis of $\mathcal I(\widehat \X)$ which corresponds to the admissible perturbation  $\widehat\X=\{(0.98{\overline3},1),\;(3.0{\overline 3},2),\; (5.08{\overline 3},3)\}$  consisting of aligned points, we conclude that the points $\X$ are misaligned  because of data inaccuracy.\hfill$\diamondsuit$}
\end{ex}

\medskip
In Example~\ref{ex2}  a set $\X$  of $20$ points close to  a  circumference is processed and the NBM algorithm detects this geometrical configuration. 

\begin{ex}\label{ex2}
{\rm 
Let  $\X$  be a set of points obtained varying the coordinates of  a set of $20$ points lying on the circumference $C$ of equation  $x^2+y^2-1=0$,  with componentwise perturbations less than $ 10^{-4}$.

The $\sigma$-Gr\"obner basis of $\mathcal I(\X)$  does not detect that the points $\X$ are close to a circumference. 
On the contrary,  the NBM algorithm, processing the set $\X$ with tolerance $\varepsilon=(0.0001, 0.0001)$, computes the stable quotient basis $\OI$
$$
\OI=\{1, y, x, y^2, xy, y^3, xy^2, y^4, xy^3, y^5, xy^4,y^6, xy^5, y^7, xy^6, y^8, xy^7, y^9, xy^8, y^{10}\} 
$$
and the set $\G$ of polynomials
\begin{eqnarray*}
\G=\left\{
\begin{array}{rcl}
g_1&=&x^2+0.99999y^2-1.00002\\
g_2&=&xy^9-2.00006xy^7+1.31256xy^5-0.31251xy^3+0.01953xy\\
g_3&=&y^{11}-3.00006y^9+3.3126y^7-1.6251y^5+0.3320y^3-0.0195y 
\end{array}
\right.
\end{eqnarray*}
The set $\G$ is not a Gr\"obner basis, but $g_2$ and $g_3$ vanish at $\X$ and $g_1$, with  coefficient vector $c_1$, is almost vanishing  at $\X$ since $\|g_1(\X)\|/\|c_1\| \approx10^{-4}$. 

Moreover, since the  coefficient vector of $g_1$  are close to those of the circumference $C$ we conclude  that the elements of $\X$ are ``almost lying" on~$C$.
\hfill$\diamondsuit$}\end{ex}

\medskip
In Example~\ref{eccezioni} the NBM algorithm processes the same set of points with two different tolerances. In the first case it detects two incompatible geometrical configurations close to  $\X$. In the second case, choosing a smaller tolerance, the NBM algorithm computes a set $\G$ very similar to a Gr\"obner basis of $\mathcal I(\X_1)$, where $\X_1$  is an admissible perturbation of $\X$.
\begin{ex}\label{eccezioni}
{\rm
Given the set $\X$ of points
\begin{eqnarray*}
{\X}=\{(1,6),(2,3),(2.449,2.449),(3,2),(6,1)\}
\end{eqnarray*}
we consider two different tolerances.

Firstly if $\varepsilon=(0.018,0.018)$ the NBM algorithm computes the stable quotient basis~$\OI=\{1, y, x, y^2, y^3\}$ and the set $\G$ of polynomial which is not a basis of a vanishing ideal since its zero set is empty:
\begin{eqnarray*}
\G=\left\{\begin{array}{rcl}
g_1&=&xy+0.00008y^2-0.00064x-0.00125y-5.99501\\
g_2&=&x^2+0.99199y^2-11.94095x-11.88550y+46.54436\\
g_3&=&y^4-14.477y^3+76.7241y^2-14.8620x-188.4194y+214.3446
\end{array}\right.
\end{eqnarray*}
In this case $g_1$ and $g_2$ highlight that the points of $\X$  almost lye on  the hyperbola $xy-6$ and on the circumference $x^2+y^2-12 x-12y +47$. In fact we have that both sets of points
\begin{eqnarray*}
\X_1&=&\{(1,6),(2,3),(\sqrt{6},\sqrt{6}),(3,2),(6,1)\} \;\;\; {\rm and}\\
\X_2&=&\{(1,6),(2,3),(6-2.5\sqrt{2},6-2.5\sqrt{2}),(3,2),(6,1)\} 
\end{eqnarray*}
are admissible perturbations of $\X$.  Nevertheless the configurations corresponding to $\X_1$ and $\X_2$ are incompatible, since $\#\X=5$ while the intersection between an hyperbola and a circumference consists of at most  $4$ points.

\smallskip
If we choose a smaller tolerance, the configuration of points near to the circumference is not detected by the algorithm. In fact, if $\varepsilon=(0.001,0.001)$ we obtain the stable quotient basis~$\OI=\{1, y, x, y^2, x^2\}$ and the set $\G$ of polynomial, with an empty zero set:
\begin{eqnarray*}
\G=\left\{\begin{array}{rcl}
g_1&=&xy+0.00008y^2-0.00064x-0.00125y-5.9950\\
g_2&=&y^3-2.3444x^2-14.3444y^2+34.1336x+75.1336y-182.1901\\
g_3&=&x^3-14.3444x^2-2.3444y^2+75.1336x+34.1336y-182.1901
\end{array}\right.
\end{eqnarray*}
In this case the $\sigma$-Gr\"obner basis $ GB_1$ of $\mathcal I(\X_1)$
\begin{eqnarray*}
GB_1=\left\{\begin{array}{l}
xy - 6\\
y^3-2.4494x^2-14.4494y^2+35.3938x+76.3938y-187.1260\\
x^3-14.4494x^2-2.4494y^2+76.3938x+35.3938y-187.1260
\end{array}\right.
\end{eqnarray*}
 consists of polynomials ``similar" to the elements of $\G$. Since $\X_1$ is an admissible perturbation also w.r.t.~the tolerance $(0.01,0.01)$ then  $\G$ highlight that the points $\X$ almost lye on a hyperbola. \hfill$\diamondsuit$ 
}
\end{ex}

\medskip
Example~\ref{ex4} shows that the term ordering $\sigma$ is only a computational tool for building a factor closed set. In fact, given the set $\X$, the NBM algorithm computes the order ideal  $\OI$ which   cannot be obtained by the exact Buchberger-M\"oller  algorithm working on $\X$ with any term ordering. 

\begin{ex}\label{ex4}
{\rm
Let $\X=\{(1.1,1.1), \;(0.9,-1.1), \;(-0.9,0.9),\; (-1.1, -0.9)\}$ be  the input points and let $\varepsilon=(0.12,0.12)$ be the tolerance. Since the vector space $P/\mathcal I(\X)$ has dimension $4$, the possible quotient bases are
\begin{eqnarray*}
&\OI_1=\{ 1,x,x^2,x^3\} \;\;\;\OI_2=\{ 1,y,y^2,y^3\} &\\
&\OI_3=\{1,y,x,x^2 \} \;\;\;\OI_4=\{1,y,x,y^2 \} \;\;\;\OI_5=\{1,y,x,xy \} &
\end{eqnarray*}
Each quotient basis  $\OI_j$, $j=1\dots 4$, is associated to the $\sigma_j$-Gr\"obner basis of $\mathcal I(\X)$, where
$\sigma_1=\sigma_2=Lex$ with  $y>x$ or $x>y$ respectively, and $\sigma_3=\sigma_4=DegLex$  with  $y>x$ or $x>y$ respectively. Nevertheless these sets are not stable quotient bases, since each evaluation matrix $M_{\OI_j}(\Xtilde)$, $j=1\dots 4$, is singular for the admissible perturbation $\Xtilde=\{(1,1), \;(1,-1), \;(-1,1),\; (-1, -1)\}$.

Vice versa $\OI_5$, computed by the NBM algorithm,  cannot be obtained using the exact Buchberger-M\"oller algorithm w.r.t.~any term ordering $\tau$. In fact the vector $t(\X)$, with $t=x^2$ or $t=y^2$ is independent of $\big\{r(\X) : r\in \{1,x,y\}\big \}$ so that  $\OI_3$, if $x^2<_\tau  xy$,  or $\OI_4$, if  $y^2<_\tau xy$, is built. 
It follows that the set $\G$ computed by the NBM algorithm 
\begin{eqnarray*}
\G=\left \{
\begin{array}{l}
y^2-0.19998x+0.01980y-1.01\\
x^2-0.20199xy+0.00201x+0.01999y-0.98980,   
\end{array}
\right.
\end{eqnarray*}
is not the $\tau$-Gr\"obner basis of $\mathcal I(\X)$, for any term ordering $\tau$. Nevertheless, $\G$ is the $\sigma$-Gr\"obner basis of $\mathcal I(\overline \X)$, where the zero set $\overline \X$ of  $\G$ is the admissible perturbation:
\begin{eqnarray*}
\;\;\;\;\;\overline \X=\{(1.099,1.099),(0.899, -1.100),(-0.899,0.901),(-1.099, -0.898)\}\;\;\;\;\;\diamondsuit
\end{eqnarray*}
 }\end{ex}

\vspace{3 mm}\noindent {\bf Acknowledgements.}
Part of this work was conducted during the Special Semester on Gr\"obner Bases supported by RICAM (the Radon Institute for Computational and Applied Mathematics, Austrian Academy of Science, Linz) and organized by RICAM and RISC (Research Institute for Symbolic Computation, Johannes Kepler University, Linz, Austria) under the scientific direction of Professor Bruno Buchberger.

The author would like to thank Prof.~H.~M.~M\"oller and Prof.~H.~Stetter for their  useful and constructive remarks.

\end{document}